\documentclass{amsart}

\usepackage{amsthm, amsfonts, amssymb, amsmath, latexsym, enumerate, times}
\usepackage[latin1]{inputenc}
\usepackage{mathrsfs}
\usepackage{mathtools}
\usepackage{letltxmacro}
\LetLtxMacro\orgvdots\vdots
\LetLtxMacro\orgddots\ddots

\makeatletter
\DeclareRobustCommand\vdots{%
	\mathpalette\@vdots{}%
}
\newcommand*{\@vdots}[2]{%
	\sbox0{$#1\cdotp\cdotp\cdotp\m@th$}%
	\sbox2{$#1.\m@th$}%
	\vbox{%
		\dimen@=\wd0 %
		\advance\dimen@ -3\ht2 %
		\kern.5\dimen@
		\dimen@=\wd2 %
		\advance\dimen@ -\ht2 %
		\dimen2=\wd0 %
		\advance\dimen2 -\dimen@
		\vbox to \dimen2{%
			\offinterlineskip
			\copy2 \vfill\copy2 \vfill\copy2 %
		}%
	}%
}
\DeclareRobustCommand\ddots{%
	\mathinner{%
		\mathpalette\@ddots{}%
		\mkern\thinmuskip
	}%
}
\newcommand*{\@ddots}[2]{%
	\sbox0{$#1\cdotp\cdotp\cdotp\m@th$}%
	\sbox2{$#1.\m@th$}%
	\vbox{%
		\dimen@=\wd0 %
		\advance\dimen@ -3\ht2 %
		\kern.5\dimen@
		\dimen@=\wd2 %
		\advance\dimen@ -\ht2 %
		\dimen2=\wd0 %
		\advance\dimen2 -\dimen@
		\vbox to \dimen2{%
			\offinterlineskip
			\hbox{$#1\mathpunct{.}\m@th$}%
			\vfill
			\hbox{$#1\mathpunct{\kern\wd2}\mathpunct{.}\m@th$}%
			\vfill
			\hbox{$#1\mathpunct{\kern\wd2}\mathpunct{\kern\wd2}\mathpunct{.}\m@th$}%
		}%
	}%
}
\makeatother

\usepackage{array}
\usepackage{comment}
\usepackage{paralist}
\usepackage{xcolor}

\newtheorem{theorem}{Theorem}
\newtheorem{lemma}[theorem]{Lemma}

\newtheorem{corollary}[theorem]{Corollary}
\newtheorem{proposition}[theorem]{Proposition}

\theoremstyle{definition}

\newtheorem{remark}[theorem]{Remark}
\newtheorem{example}[theorem]{Example}

\newcommand{\calE}{{\mathcal E}}
\newcommand{\calL}{{\mathcal L}}
\newcommand{\calM}{{\mathcal M}}
\newcommand{\calP}{{\mathcal P}}
\newcommand{\calG}{{\mathcal G}}
\newcommand{\calO}{{\mathcal O}}
\newcommand{\calX}{{\mathcal X}}
\newcommand{\calY}{{\mathcal Y}}
\newcommand{\bbP}{{\mathbb P}}
\newcommand{\bbR}{{\mathbb R}}
\newcommand{\bbZ}{{\mathbb Z}}

\newcommand{\Pic}{\operatorname{Pic}}
\newcommand{\Uncoll}{\operatorname{Uncoll}}
\def\geq{\geqslant}
\def\leq{\leqslant}
\def\le{\leqslant}
\def\ge{\geqslant}

\begin{document}
 
\title[Irrational rays]{Irrational nef rays at the boundary of the Mori cone for very general blowups of the plane}

\author{Ciro Ciliberto}
\address{Dipartimento di Matematica, Universit\`a di Roma Tor Vergata, Via O. Raimondo
 00173 Roma, Italia}
\email{cilibert@axp.mat.uniroma2.it}

\author{Rick Miranda}
\address{Department of Mathematics, Colorado State University, Fort Collins (CO), 80523,USA}
\email{rick.miranda@colostate.edu}
 
\author{Joaquim Ro\'e}
\address{Departament de Matem\`atiques, Universitat Aut\`onoma de Barcelona, 08193 Bellaterra (Barcelona), Catalunya}
\email{jroe@mat.uab.cat}

 
 
\maketitle

\begin{abstract} In this paper we develop a technique for discovering (non-effective) irrational rays at the boundary of the Mori cone for linear systems on a general blowup of the plane, and give examples of such irrational rays.
\end{abstract}

\section*{Introduction} 

Let $X_s$ be the blow-up of the complex projective plane at $s$ very general points. 
We will usually assume $s\geq 10$.
Let $\calL=\calL_d(m_1,\ldots, m_s)$, with $d>0$,
be the linear system on $X_s$
corresponding to plane curves of degree $d$
having multiplicities at least $m_1,\ldots,m_s$ at the given points
(we will use exponential notation for repeated multiplicities).
If we denote by $H$ the class of the pullback to $X_s$ of a general line in the plane,
and by $E_i$ the exceptional divisor over the $i$-th point that is blown up,
then the linear system $\calL$ corresponds to
the global sections of the line bundle $\calO_{X_s}(dH-\sum_i m_i E_i)$.
The (projective) dimension of the linear system is
$\dim (H^0(X_s, \calO_{X_s}(dH-\sum_i m_i E_i)))-1$.
The system is said to be \emph{effective} if the dimension is non-negative,
i.e., there are effective divisors in the linear system.
If the dimension is $-1$, we will say that the system is \emph{non-effective},
or \emph{empty}.

We note that the Picard group of $X_s$
is the free abelian group of rank $s+1$ 
generated by the classes of $H$ and $E_1,\ldots, E_s$.
We define $N=N^1(X_s)$ to be $\Pic(X_s)\otimes_{\bbZ}\bbR$;
this is a real vector space of dimension $s+1$.

A \emph{ray} in $N$ is the set of all non-negative real multiples of a non-zero vector in $N$.
Many of the concepts applicable to elements of the Picard group $\Pic(X_s)$
may be extended to $N$.
A ray is called \emph{effective} if there is a (necessarily integral) vector in the ray
that represents an effective divisor class in $\Pic(X_s)$.
The ray is called \emph{rational} if it contains a nonzero integral vector
(i.e., an element of the Picard group).
The \emph{degree} $\deg(R)$ of a ray $R$,
the coefficient of $H$,
is not well-defined,
but its sign is;
hence we may speak of a ray of positive/zero/negative degree.
Every effective ray $R$ must have non-negative degree.
Similarly, if $R$ and $R'$ are two rays, the intersection $R\cdot R'$ is not well defined,
but its sign is. In particular this applies to the
self-intersection of a ray.  
Any rational ray with $\deg(R) > 0$ and $R^2 > 0$ is effective
(by Riemann-Roch).

The \emph{effective cone} is the cone generated by effective rays,
i.e., all finite linear combinations of effective divisors with positive real coefficients.
In general the effective cone is not closed;
its closure is called the \emph{Mori cone}.
The dual of the Mori cone is the cone of \emph{nef} divisors, or nef rays;
these are divisors/rays that intersect all effective divisors non-negatively.
This nef cone is closed.

In \cite{CHMR13}, a \emph{good ray} was defined
as a non-effective rational ray with non-negative degree (intersection with $H$)
and zero self-intersection.
Such a ray is extremal for the Mori cone and the nef cone, and is nef
(see \cite[Lemma 3.8]{CHMR13}).

In that paper we also defined a \emph{wonderful} ray,
as a nef ray that has self-intersection zero and is irrational.
Up to this point, no wonderful ray (with all coordinates nonzero)
for any $s \geq 10$ has been detected.
Some of the motivation for discovering wonderful rays is as follows.

Define the \emph{De Fernex ray} $F_s$
to be the ray generated by $\sqrt{s-1}H-\sum_{i=1}^s E_i$.
A ray $R$ is said to be \emph{De Fernex positive, negative or orthogonal} 
according to $R\cdot F_s$ being positive, negative or null.
The \emph{Strong $\Delta$-Conjecture}
(see \cite[Conjecture 3.10]{CHMR13})
is that if $s \geq 11$,
and $R$ is a rational De Fernex non--positive ray of self-intersection zero,
then $R$ is not effective,
and therefore is a good ray.
(See \cite{DeF01}; there is a refinement for the $s=10$ case.)
Note that $R\cdot F_s \leq 0$ implies $R\cdot K_s>0$,
where $K_s = -3H+\sum_i E_i$ is the canonical divisor on $X_s$.

The Strong $\Delta$-Conjecture would imply that
an irrational De Fernex non--positive ray $R$ with self-intersection zero
is a wonderful ray,
since any such ray would then be a limit of good rays which are nef.

This remark implies, in particular, that the \emph{Nagata ray} $\sqrt{s}H-\sum_i E_i$
would be wonderful (if $s$ is not a square), 
and this would prove the Nagata Conjecture (see \cite{N59}).
Of course the Nagata Conjecture is very important 
and has been shown for decades that it is difficult to prove.
However, in view of the Strong $\Delta$-Conjecture,
from a conceptual Mori cone viewpoint,
proving that the Nagata ray is wonderful
is not philosophically more important
than proving that any other De Fernex negative irrational ray
with self-intersection zero is wonderful.

As noted above, no wonderful ray (with all coordinates nonzero)
for any $s \geq 10$ has been detected up to now.
In this paper we fill this gap by proving the following.

\begin{theorem}\label{thm:main}
For all  $s \geq 10$, wonderful rays $R$ exist. 
For all  $s \geq 13$, wonderful rays $R$ with $R\cdot K_s>0$ exist. 
For $s=14$, for all $s\ge 13$ such that $s-4$ is a square, 
and for all $s\ge 18$ such that $s-2$ is a square,
De Fernex negative wonderful rays exist.
\end{theorem}

Our proof is explicit, in that for each $s$, 
we exhibit an irrational ray with self-intersection zero, 
all of whose coordinates are positive, 
and prove that it is wonderful by proving that it is a limit of explicit good rays 
(note that all good rays known so far were isolated, 
so this is also the first example of accumulation of selfintersection zero classes 
on the boundary of the Mori cone).
For certain values of $s$, as stated in Theorem \ref{thm:main},
these wonderful rays intersect the canonical divisor positively,
and intersect the De Fernex ray negatively.

Specifically, we start by exploiting an infinite sequence of Cremona transformations,
which we then apply to a carefully chosen good ray,
and show that the limit ray exists and is irrational.
Since Cremona transformations preserve the 'goodness' of the rays,
the limit is therefore a limit of good rays, and is therefore wonderful.
The rays obtained in this way are orthogonal to $K_s$.

To produce wonderful rays that meet $K_s$ positively,
and are De Fernex-negative,
we exploit a degeneration technique
(described in Section \ref{sec:collision})
that allows us to coherently reduce $K_s$-positive and De Fernex-negative systems
to the systems generating the sequence of good rays found in the first step above.
The limit of these systems then provide the examples
which prove Theorem \ref{thm:main}.

The construction that we present and develop 
is a general technique for iteratively generating good rays, 
based on a judicious use of
degeneration procedures and Cremona transformations.
This technique can be applied more extensively,
to produce additional wonderful rays, 
some of which have additional properties;
we are developing this.
In particular using this same strategy, 
we can get many other irrational points on the boundary of the Mori cone,
though we do not present them in this paper.

The linear systems we consider 
will have points with at most three distinct multiplicities,
i.e., they will be of the form $\calL_d(a^{s_a},b^{s_b},c^{s_c})$. 
This is only a technical device to simplify computations;
clearly many other examples can be found with more different multiplicities.

There is a close relationship between our results here
and the possible existence of irrational Seshadri constants.
Namely, there are many rays of selfintersection zero
(like the Nagata ray)
which, if proven to be wonderful,
would imply that some Seshadri constant 
at a blow-up of $\bbP^2$ is irrational.
The wonderful rays presented here are not among these; 
however, it is possible that our techniques can be used 
to show the existence of such rays.
(See 
\cite{HR08}, \cite{CHMR13}, \cite{DKMS16} and \cite{HH18} 
for more details.)

In the first two sections, we present the technical tools used to obtain wonderful rays. In the last three, we prove separately the three existence claims of our main Theorem.

\medskip

{\bf Acknowledgements:} Ciro Ciliberto is a member of GNSAGA of the Istituto Nazionale di Alta Matematica ``F. Severi".
We thank the Centre de Recerca Matem\`atica 
for the hospitality in arranging Ciliberto's and Miranda's visit to Barcelona within its ``Research in Pairs'' program. 
This work was also partially supported by project PID202-116542GB-I00 from the Spanish Ministerio de Ciencia e Innovaci\'on.

\section{Collision of $r^2$ points}
\label{sec:collision}

A key part of the construction 
involves a degeneration of a linear system $\calL_d(m_1,\ldots, m_s)$
where $r^2$ of the $s$ points, of equal multiplicity $m$, come together;
following \'Evain, who treated particular cases of the situation in the 90's \cite{E98}, 
we call this a \emph{collision}.
We will index the points so that the first $r^2$ multiplicities are all equal to $m$.
The analysis of this situation for $r=2$ was developed in \cite{CM05};
here we need to allow any integer $r\ge 2$, 
but the same technique 
combined with Nagata's result of \cite{N59} on square numbers of points 
is enough for our purposes. 

We consider a trivial family $\calX = X_{s-r^2} \times \Delta$ over a disc $\Delta$,
and blow up a general point in the central fiber over $0\in \Delta$ 
to obtain the threefold $\calX'$.
This produces a degeneration of $X_{s-r^2}$
to a union of two surfaces, 
a plane (the exceptional divisor for the blowup)
and the proper transform $F$ of the original $X_{s-r^2}$ fiber,
which is now isomorphic to $X_{s-r^2+1}$.
These two surfaces intersect transversely along a smooth rational curve $R$
which is a line in the plane
and a $(-1)$-curve in $F$.

We now choose $r^2$ general points on the plane;
extend these $r^2$ general points to the general fiber 
using $r^2$ sections of the projection of $\calX'$  to $\Delta$,
and blow up those $r^2$ sections to ruled surfaces
$\mathcal E_1,\ldots,\mathcal E_{r^2}$.
This then produces a threefold $\calY$
which is a degeneration of $X_s$,
to a union of a surface $P \cong X_{r^2}$ and $F \cong X_{s-r^2+1}$,
intersecting transversely along the double curve $R$.
This smooth rational curve $R$
is the pullback of a general line in the surface $P$
and remains a $(-1)$-curve in the surface $F$.

We have the line bundle
corresponding to $\calL_d(m_{r^2+1},\ldots,m_s)$ on $X_{s-r^2}$,
and can extend it trivially to $\calX$.
If we pull that back to the first blowup $\calX'$,
we see that this restricts to the bundle
corresponding to $\calL_d(0,m_{r^2+1},\ldots,m_s)$ on the surface $F$,
and to the trivial bundle on the plane.
We then pull that back to the second blowup $\calY$,
and tensor by $\calO_{\calY}(-tP-m \sum_{i=1}^{r^2} \calE_i)$, 
with $t$ a non--negative integer (called the \emph{twisting parameter}).
This produces a line bundle $\calM$ on $\calY$,
which restricts to the general fiber
in a bundle whose associated linear system is
the original system $\calL_d(m_1,\ldots,m_s)$.

The restriction of $\calM$ to $P$
is a bundle with associated linear system
$\calL_t(m^{r^2})$.
The restriction of $\calM$ to $F$
is a bundle with associated linear system $\calL_d(t,m_{r^2+1},\ldots,m_s)$.

At this point we choose the twisting parameter $t$
to be the minimum $t$ 
such that the linear system $\calL_t(m^{r^2})$ is effective.
By Nagata's theorem, if $r>3$ then $t>rm$, and it is not hard to see that if $r\le 3$ then $t=rm$. 
The principle of semicontinuity guarantees that
the dimension of the general linear system $\calL_d(m_1,\dots,m_s)$
is at most equal to the dimension of the linear system on the reducible surface $P+F$.
The linear system on the reducible surface
corresponds to the fiber product
of the space of sections on $P$ with the space of sections on $F$,
fibered over the restriction to the space of sections on $R$.

Let us investigate this in more detail in the cases $r = 2,3$;
we do not need a precise description if $r>3$.
We note that since we are taking the minimum twist parameter $t$
to make the system on $P$ effective,
the restriction map from the space of sections on $P$
to the space of sections on the double curve $R$ 
(which is $H^0(\calO_R(t))$)
is injective, onto a subspace $W \subset H^0(R, \calO_R(t))$ of codimension $c$.
Hence the fiber product is isomorphic to
the subspace of $H_F=H^0(F,\calO_F(dH-tR-\sum_{i=r^2+1}^s m_i E_i)$
which restricts to elements of $W$.
By the generality of the points chosen on $P$,
the subspace $W$ is transversal to the full restriction,
(see \cite[Section 3]{CM98})
and so the desired subspace also has codimension $c$ in $H_F$.


If $r=2$, the linear system on $P$ is $\calL_{2m}(m^4)$,
whose elements are sums of $m$ conics in the system $\calL_2(1^4)$.
This system restricts to $R$ in a linear system on $R$ which is not complete
(being of degree $2m$ and dimension $m$);
it has codimension $m$.
This imposes the additional matching condition on the linear system on $F=X_{s-3}$:
it must restrict to divisors on $R$ which are members of the linear system coming from $P$.
We see then that the fiber product has codimension $m$
in the system on $F$;
these are the additional 'matching' conditions to be imposed on the system on $F$
(over and above the point of multiplicity $t=2m$).

If $r=3$, the linear system on $P$ is $\calL_{3m}(m^9)$, 
which has a unique effective element,
namely the unique cubic through the $9$ general points, taken with multiplicity $m$.
This cubic meets the double curve $R$ in three (general) points,
and will impose on the other component $X_{s-8}$ three points of certain multiplicity
on the double curve.
We won't need to explicitly determine the multiplicity of these additional points;
it will be enough for our purposes to observe that the restricted system on $R$ has degree $3m$ and dimension $0$,
which means that the restriction subpace $W$ above has codimension $3m$.
Hence the number of additional 'matching' conditions imposed on $F=X_{s-8}$
is at least $3m$ (over and above the point of multiplicity $t=3m$).

In particular,
whenever the linear system on $F$ (with the matching conditions) is empty,
the collision shows that the original system was empty.
It will be useful for us to proceed in reverse:
from a system on the plane with $s-r^2+1$ points,
whose dimension is known,
to another system with $s$ points obtained replacing the first point,
of multiplicity $rm$, by $r^2$ points of multiplicity $m$. 
We call this an \emph{uncollision step}.

We summarize this in the following.

\begin{lemma}\label{lem:uncollision}
Fix $r=2$ or $3$, $s \geq r^2+1$,
and multiplicities $m$, $m_{r^2+1},\ldots,m_s$.
\begin{itemize}
\item[(a)] If $r=2$ and $\dim \calL_d(2m,m_{5},\ldots,m_s) < m$,
then $\calL_d(m^4,m_{5},\ldots,m_s)$ is empty.
\item[(b)] If $r=3$ and $\dim \calL_d(3m,m_{10},\ldots,m_s) < 3m$,
then $\calL_d(m^9,m_{10},\ldots,m_s)$ is empty.
\item[(c)] If $r \geq 4$ and $\calL_d(rm+1,m_{r^2+1},\ldots,m_s)$ is empty,
then the uncollided system $\calL_d(m^{r^2},m_{r^2+1},\ldots,m_s)$ is also empty.
\end{itemize}
\end{lemma}

Three comments are in order.
First, the above reductions are sharp for $r=2$,
but not for $r=3$ or $r \geq 4$:
there are additional matching conditions in order for a curve
in the degenerate surface $P+F$
to be a limit of a curve on the general fiber.

Second, we note that for both $r=2,3$,
the systems before and after the collision
have the same self-intersection.
In particular, if one is zero, so is the other;
this will be important in our application.

Third, the process of considering an 'uncollision' behaves well with taking limits of rays.
Given a linear system $\calL$, and an index $i$ denoting one of the multiplicities,
we may define the \emph{uncollision} $\Uncoll_r(\calL,i)$
as the system replacing the $i$-th multiplicity $m_i$
by $r^2$ points of multiplicity $m_i/r$.
This makes sense at the level of linear systems if $m_i$ is divisible by $r$,
but also makes sense as elements of $N$,
and additionally makes sense for rays in $N$.
In particular, if $\calL_k$ is a sequence of linear systems,
giving rise to rays $[\calL_k] \in N$,
then
\begin{equation}\label{uncoll-limit}
\lim_{k\to\infty} [\Uncoll_r(\calL_k,i)] = [\Uncoll_r(\lim_{k \to\infty}(\calL_k),i)]
\end{equation}
as rays in $N$.  We finally note that the uncollision process is given by rational parameters, and so preserves rationality and irrationality of rays.

\section{Useful Cremona transformations}
\label{sec:Cremonas}

In this section we explore some Cremona maps 
that act on linear systems of the form $\calL_d(a^{s_a},b^{s_b},c^{s_c})$.

Consider the Cremona--Kantor (CK) group $\calG_s$
generated by quadratic transformations 
based at $n$ general points $x_1 , \dots, x_s$ of the plane
and by permutations of these points
(see \cite{DuV36}, \cite{A02}, \cite[Chapter 7]{Dol12}).
The group $\calG_s$ acts on the set of linear systems
of the type $\calL_{d}(m_1, \dots, m_s )$.
All systems in the same (CK)-orbit (or (CK)-equivalent) have
the same expected, virtual and true dimension.
A linear system $\calL_{d}(m_1, \dots, m_s )$ is \emph{Cremona reduced}
if it has minimal degree in its (CK)-orbit,
and this is the case if and only if the degree 
is greater or equal to the sum of the three largest multiplicities
(see \cite[p. 402-402, Thms 8 and 10]{C31}).

An element $\phi \in \calG_s$,
seen as a linear automorphism of 
$\Pic(X_s)=\bbZ H\oplus \bbZ E_1\oplus \dots\oplus E_s$,
can be specified by giving its \emph{characteristic matrix},
i.e., the matrix with respect to the standard basis $(H, E_1,\dots, E_s)$. 
The \emph{homaloidal net} of $\phi$
(i.e., the pullback of the net of lines by $\phi$) is $\calL_{d}(m_1, \dots, m_s)$
where $(d, -m_1, \dots, -m_s)$ is the first column of the characteristic matrix of $\phi$.

\begin{example}
There exist four homaloidal types with homogeneous multiplicities
(see \cite[2.5.5]{A02}, \cite[7.2.2]{Dol12}).
The simplest one corresponds to the \emph{(standard) quadratic Cremona map}
on three points, whose characteristic matrix is 
\[
	Q=
	\left(\begin{array}{rrrr}
	2 & 1 & 1 & 1 \\
	-1 & 0 & -1 & -1 \\
	-1 & -1 & 0 & -1 \\
	-1 & -1 & -1 & 0
	\end{array}\right).
\]
Any permutation of the $3$ rightmost columns 
gives rise to a \emph{distinct} quadratic Cremona map; 
the given matrix corresponds to the involutive quadratic map.
	
Of course, the quadratic Cremona map can be applied 
on any subset of three points among the set of points $\{p_1,\dots,p_s\}$
 that we blow up, for any $s\ge 3$;
the characteristic matrix of the corresponding element in $\calG_s$
is obtained from $Q$ by adding suitably many rows and columns of the identity matrix. 
Similarly, any Cremona map defined for $X_s$ 
can be applied to $X_{s'}$ with $s'>s$ 
by selecting a suitable set of $s$ points among the $s'$, 
and the matrix is obtained by adding rows and columns of the identity.
	
The other three homogeneous homaloidal types 
are attributed to Sturm, Geiser and Bertini; 
as in the quadratic case, 
for each of them there is a unique $\phi\in\calG_s$ 
with that type and of order $2$.
Their characteristic matrices are
\[
S=
	\begin{psmallmatrix*}[r]
	5 & 2 & 2 & 2 & 2 & 2 & 2 \\
	-2 & 0 & -1 & -1 & -1 & -1 & -1 \\
	-2 & -1 & 0 & -1 & -1 & -1 & -1 \\
	-2 & -1 & -1 & 0 & -1 & -1 & -1 \\
	-2 & -1 & -1 & -1 & 0 & -1 & -1 \\
	-2 & -1 & -1 & -1 & -1 & 0 & -1 \\
	-2 & -1 & -1 & -1 & -1 & -1 & 0
	\end{psmallmatrix*},
G=
	\begin{psmallmatrix*}[r]
	8 & 3 & 3 & 3 & 3 & 3 & 3 & 3 \\
	-3 & -2 & -1 & -1 & -1 & -1 & -1 & -1 \\
	-3 & -1 & -2 & -1 & -1 & -1 & -1 & -1 \\
	-3 & -1 & -1 & -2 & -1 & -1 & -1 & -1 \\
	-3 & -1 & -1 & -1 & -2 & -1 & -1 & -1 \\
	-3 & -1 & -1 & -1 & -1 & -2 & -1 & -1 \\
	-3 & -1 & -1 & -1 & -1 & -1 & -2 & -1 \\
	-3 & -1 & -1 & -1 & -1 & -1 & -1 & -2
	\end{psmallmatrix*}
\]
and
\[
B=
	\begin{psmallmatrix*}[r]
	17 & 6 & 6 & 6 & 6 & 6 & 6 & 6 & 6 \\
	-6 & -3 & -2 & -2 & -2 & -2 & -2 & -2 & -2 \\
	-6 & -2 & -3 & -2 & -2 & -2 & -2 & -2 & -2 \\
	-6 & -2 & -2 & -3 & -2 & -2 & -2 & -2 & -2 \\
	-6 & -2 & -2 & -2 & -3 & -2 & -2 & -2 & -2 \\
	-6 & -2 & -2 & -2 & -2 & -3 & -2 & -2 & -2 \\
	-6 & -2 & -2 & -2 & -2 & -2 & -3 & -2 & -2 \\
	-6 & -2 & -2 & -2 & -2 & -2 & -2 & -3 & -2 \\
	-6 & -2 & -2 & -2 & -2 & -2 & -2 & -2 & -3
	\end{psmallmatrix*}
\]
	respectively. 
\end{example}

\begin{example}
In addition to the homogeneous Cremona maps of the previous example,
we shall use as building blocks two families of \emph{quasi-homogeneous} involutions. 
The first is the \emph{de Jonqui\`eres map} on $2n+1$ points with characteristic matrix
\[
J_n=
	\begin{pmatrix*}[r]
	1+n & n & 1 & 1 & \cdots & 1 \\
	-n & 1-n & -1 & -1 & \cdots & -1 \\
	-1 & -1 & -1 & 0 & \cdots & 0 \\
	-1 & -1 & 0 & -1 & \ddots & \vdots \\
	\vdots & \vdots & \vdots & \ddots & \ddots & 0 \\
	-1 & -1 & 0 & \cdots & 0 & -1\end{pmatrix*}
\]
(this matrix is not explicitly given in \cite[2.6.10, 3.4.3]{A02} or \cite[7.2.3]{Dol12},
but it is easy to recover it from the fact that the de Jonqui\`eres Cremona map 
is the composition of $n$ quadratic Cremona transformations
based at points $\{p_1,p_{2i},p_{2i+1}\}$, for $i=1,\ldots,n$,
and permuting some points).
	
The second family have an even number of base points
and does not seem to have received any special attention so far. 
The characteristic matrix defining the map with $2n+2$ base points is
\[
C_n=
	\begin{pmatrix*}[r]
	1+n^2 & -n+n^2 & n & n & \cdots & n \\
	n-n^2 & 2n-n^2 & 1-n & 1-n & \cdots & 1-n \\
	-n & 1-n & 0 & -1 & \cdots & -1 \\
	-n & 1-n & -1 & 0 & \ddots & \vdots \\
	\vdots & \vdots & \vdots & \ddots & \ddots & -1 \\
	-n & 1-n & -1 & \cdots & -1 & 0
	\end{pmatrix*}.
\]
The reader may check that this is the result of applying a de Jonqui\`eres map
based at points $p_2, p_3, \dots, p_{2n+2}$
followed by another de Jonqui\`eres
based at $p_1, p_3,\dots, p_{2n+2}$
(where $p_3,\dots,p_{2n+2}$ are the simple points of each de Jonqui\`eres homaloidal net).
\end{example}

\begin{example}\label{exa:iterable_cremonas}
	We will use in our constructions some particular familes of Cremona maps of infinite order, built composing maps of the previous kinds.
	The maps in the first family involve $s=2n+7$ points; they are obtained by composing a quintic Sturm map based on points $p_{2n+2}, \dots, p_{2n+7}$ followed by a de Jonqui\`eres map based on (disjoint) points $p_1, \dots, p_{2n+1}$. 
	Extending the matrices $J_n$ and $S$ with the suitable number of rows and columns and multiplying, one obtains the corresponding characteristic matrix, namely
\[
\mathit{JS}_n=
\begin{psmallmatrix*}[r]
5+5n & n & 1 & 1 & \cdots & 1 & 2+2n& 2+2n& 2+2n& 2+2n& 2+2n& 2+2n\\
-5n & 1-n & -1 & -1 & \cdots & -1 & -2n& -2n& -2n& -2n& -2n& -2n\\
-5 & -1 & -1 & 0 & \cdots & 0 &-2&-2&-2&-2&-2&-2\\
-5 & -1 & 0 & -1 & \ddots & \vdots& \vdots& \vdots& \vdots& \vdots& \vdots& \vdots \\
\vdots & \vdots & \vdots & \ddots & \ddots & 0 & \vdots& \vdots& \vdots& \vdots& \vdots& \vdots\\
-5 & -1 & 0 & \cdots & 0 & -1&-2&-2&-2&-2&-2&-2\\
-2 & 0 & 0 & \cdots & 0 & 0 & 0 & -1 & -1 & -1 & -1 & -1 \\
-2 & 0 & 0 & \cdots & 0 & 0 & -1 & 0 & -1 & -1 & -1 & -1 \\
-2 & 0 & 0 & \cdots & 0 & 0 & -1 & -1 & 0 & -1 & -1 & -1 \\
-2& 0 & 0 & \cdots & 0 & 0  & -1 & -1 & -1 & 0 & -1 & -1 \\
-2 & 0 & 0 & \cdots & 0 & 0 & -1 & -1 & -1 & -1 & 0 & -1 \\
-2 & 0 & 0 & \cdots & 0 & 0 & -1 & -1 & -1 & -1 & -1 & 0
\end{psmallmatrix*}.
\]
In particular, the homaloidal systems of these maps 
have the form $\calL_{5+5n}(5n,5^{2n},2^6)$.
We note that if we transform a linear system $\calL$
whose parameters have the 'shape' $\calL_d(a,b^{2n},c^6)$,
then the result is a linear system with that same shape.

An additional useful family, 
involving $s=2n+8$ points for $n\ge 1$, 
can be obtained by composing an octic Geiser map 
based on points $p_{2}, \dots, p_{8}$ 
followed by a quasi-homogeneous map with characteristic matrix $C_{n+3}$
based on all points $p_1, \dots, p_{2n+8}$. 
Extending the matrix $S$ with the suitable number of rows and columns 
and multiplying, 
one obtains the corresponding characteristic matrix $CG_n$,
which we omit for brevity,
noting only that its homaloidal system is 
$\calL_{8 n^{2} + 27 n + 17} \left(8 n^{2} + 19 n + 6,(8 n + 6)^7,(8 n + 3)^{2n}\right)$.

\end{example}

%


\section{Wonderful rays in $K_s^\perp$}\label{sec:wonderfulKperp}

Because the intersection of the Mori cone 
with any hyperplane of the form $d= {}$constant is compact, 
any infinite set of rays on the Mori cone has some accumulation ray. 
This simple observation is enough to provide 
many interesting rays on the boundary of the Mori cone.
Indeed, if $s\ge 9$ then every divisor class not multiple of $K_s$ 
has an infinite orbit under the action of the Cremona group, 
and if $R$ is a ray of selfintersection zero on the boundary of the Mori cone, 
then every ray in its orbit also has selfintersection zero 
and lies on the boundary of the Mori cone. 
As a consequence there exist accumulation rays of selfintersection zero 
on the boundary of the Mori cone, 
an important fact which has not been observed before.
A careful choice of Cremona maps allows us 
to obtain explicit \emph{irrational} rays with such properties:

\begin{proposition}\label{pro:wonderful-cremona}
Let $n$ be an integer, and consider
$\alpha_n=\sqrt{n(n - 1)}$ and $\beta_n=\sqrt{49n^2-28}$.
The rays generated by 
$$
W_{\textit{odd}} = (5 n^{2} + 4 n) H
- n\left(3 n + 2  \alpha_n \right)E_1
- \left( 3 n + 2  \alpha_n\right) \sum_{i=2}^{1+2n}E_i
- n\left(2+n-\alpha_n\right)\sum_{i=2n+2}^{2n+7}E_i
$$ 
on $X_{7+2n}$ if $n\ge 2$ and by 

\begin{gather*}
W_{\textit{even}} =
14n \, {\left(8  n^{2} + 27  n + 16\right)} H
- 7n \left(n + 2\right)\left(9 n + \beta_n + 6 \right) E_1 \\
- n (21 n^{2} - 3 n\beta_n + 126 n - 2 \beta_n + 84)\sum_{i=2}^8 E_i 
- 7n\left(9 n + \beta_n + 6 \right)\sum_{i=9}^{2n+8} E_i
\end{gather*}
on $X_{8+2n}$ if $n\ge 1$ are wonderful.
\end{proposition}

\begin{corollary}\label{sgeq10wonderful}
	For every $s\ge 10$ there exist wonderful rays on $X_s$.
\end{corollary}

\begin{remark}
	Our methods do provide wonderful rays 
simpler than $W_\textit{even}$ on $X_s$ for even $s$. 
We choose this divisor because it will be useful later on, 
to construct De Fernex negative wonderful rays.
\end{remark}

\begin{proof}
Let us tackle the odd case first.
We are assuming that $n\ge 2$, so $\alpha_n$ is irrational;
moreover the selfintersection of the given system is readily computed to be zero.
We only have to show that it is nef, or equivalently, 
that the ray it spans is a limit of nef rays.
	
Let $\phi$ be the Cremona map with characteristic matrix $\textit{JS}_n$ 
given in Example \ref{exa:iterable_cremonas}.
The image under $\phi$ of a linear system of the form $\calL_{d}(a,b^{2n},c^{6})$
is again of that form,
and is computed multiplying the matrix $\textit{JS}_n$ by $(d,-a,(-b)^{2n},(-c)^{6})$;
the result has the parameters $(d',-a',(-b')^{2n},(-c')^{6})$ where
\[
	\begin{pmatrix}
	d'\\a'\\b'\\c'
	\end{pmatrix}=
	\begin{pmatrix}
	5 n + 5  &-n&  -2 n & -12 n - 12\\
	5 n & 1-n& -2n& -12 n  \\
	5 & -1 &-1 &-12 \\
	2 & 0 & 0 & -5  
	\end{pmatrix}
	\begin{pmatrix}
	d\\a\\b\\c
	\end{pmatrix}.
\]
So $\calL_{d}(a,b^{2n},c^{6})$ is mapped to $\calL_{d'}(a',(b')^{2n},(c')^6)$.
This matrix diagonalizes with eigenvalues $1$ and $2n\pm\alpha_n-1$;
thus applied iteration of $\phi$ to a general $(d,a,b,c)$ 
converges to the eigenspace of the dominant eigenvalue  $2n+\alpha_n-1$,
and this eigenspace is spanned by 
$$
\left(5 n^{2} + 4 n,3 n^{2} + 2  n \alpha_n,
	3 n + 2 \alpha_n, n\left(2+n-\alpha_n\right) \right).
$$
Since the vector $(1,0,0,0)$ corresponding to the ample class $\calL_1(0^{2n+7})$
does not belong to the span of the eigenspaces of eigenvalues $1$ and $2n-\alpha_n-1$,
the iterated application of $ \phi$ to $\calL_1(0^{2n+7})$
(which gives the homaloidal systems of the powers of $\phi$, obviously nef)
converges to the claimed ray.
	
	The even case is treated similarly; it is enough to show a sequence of nef classes converging to the given ray, and these nef classes will be the homaloidal classes of the powers of a suitable Cremona map. 
	In this case we use the Cremona map $\psi$ with characteristic matrix $CG_n$ of Example \ref{exa:iterable_cremonas}.
	The image under $\psi$ of a linear system of the form $\calL_{d}(a,b^{7},c^{2n})$ is  $\calL_{d'}(a',(b')^{7},(c')^{2n})$, where in this case
	\[
	\begin{pmatrix}
	d'\\a'\\b'\\c'
	\end{pmatrix}=
	\begin{pmatrix}
8 n^{2} + 27 n + 17 & -n^{2} - 5 n - 6 & -21 n^{2} - 70 n - 42 & -2 n^{2} - 6 n \\
8 n^{2} + 19 n + 6 & -n^{2} - 4 n - 3 & -21 n^{2} - 49 n - 14 & -2 n^{2} - 4 n \\
8 n + 6 & -n - 2 & -21 n - 15 & -2 n \\
8 n + 3 & -n - 2 & -21 n - 7 & -2 n + 1	\end{pmatrix}
	\begin{pmatrix}
	d\\a\\b\\c
	\end{pmatrix}.
	\]
	This matrix diagonalizes with dominant eigenvalue  $(n\beta_n  + 7 n^{2} - 2)/2$, and associated eigenspace spanned by 
	\begin{gather*}
	\left( 14n (8  n^{2} + 27  n + 16), 
	7n (n + 2)(9 n + \beta_n + 6), \right.\\
	\left. n (21 n^{2} - 3 n\beta_n + 126 n - 2 \beta_n + 84), 
	7n(9 n +  \beta_n + 6 )\right).
	\qedhere
	\end{gather*}
	\end{proof}

The reader may check that the wonderful classes of Proposition \ref{pro:wonderful-cremona} are orthogonal to the canonical divisor.
In fact, this will be the case for every wonderful ray constructed by iterating Cremona maps, because on one hand, a converging sequence of such classes necessarily have increasing degrees, and on the other hand, Cremona maps preserve the canonical class and the intersection product. This implies that $\lim (K_s\cdot(\phi^k(\calL))/\deg(\phi^k(\calL)))=0$, or, in other words, the limit ray is orthogonal to the canonical class.

\section{Wonderful rays in $K_s^+$}

By the Cone Theorem, the shape of the Mori cone on the half-space $K_s^-$ of classes which intersect the canonical divisor negatively is governed by the rays generated by $(-1)$-curves. 
On the orthogonal $K_s^{\perp}$ hyperplane this is no longer quite the case, 
as shown by the existence of wonderful rays, 
but these wonderful rays are very particular, as they are also limits of $(-1)$--rays (indeed, if $C$ is a $(-1)$--curve whose class does not belong to a certain linear space, then the class of $\phi^k(C)$ converges to the wonderful ray of Proposition \ref{pro:wonderful-cremona} as well). 
So it will be much more compelling evidence in favor of the Strong $\Delta$-Conjecture to show wonderful rays on $K_s^+$. 

To construct a sequence of good rays converging to a wonderful ray in $K_s^+$,
we will use the Cremona maps $\phi$ and $\psi$ above,
and the uncollision described in Section \ref{sec:collision}.
We are guided by the commutativity of taking limits and uncollisions,
as noted in \eqref{uncoll-limit}.

\color{black}
\begin{remark}
\label{rem:uncollision_K}
If $W$ is an $\bbR$-divisor class with $W \cdot K_s=0$ and $W'$ is obtained from $W$ by uncolliding a point of multiplicity $rm>0$ to $r^2\ge 4$ points of multiplicity $m$, then $W' \cdot K_s >0$.
Indeed, writing $W=dL-\sum m_i E_i$ we have $W\cdot K_s=\sum m_i -3d$ and
\[
W'\cdot K_s = \sum m_i - rm + r^2m -3d = W \cdot K_s + (r^2-r)m>W\cdot K_s.
\]  
\end{remark}

It would be convenient if we could simply uncollide the wonderful rays
found in Section \ref{sec:wonderfulKperp} and prove that those rays were also wonderful;
by the above remark, they would lie in $K_s^+$. 
However the collision/uncollision analysis and construction
is only available for actual linear systems,
and not for irrational rays.

Hence we must finesse this, by uncolliding each (integral) linear system in the sequence,
and show that the limit of these uncollided systems is wonderful.
Indeed, since the constructions of the uncollision are relatively simple linear transformations of the parameters, it is elementary that the limit ray of the uncollided systems
will be the formal uncollision of the wonderful ray found earlier,
as noted in \eqref{uncoll-limit}.

We can codify this approach with the following:

\begin{lemma}\label{lem:uncollwonderful}
Suppose that $\{\calL_k\}$ is a sequence of linear system rays such that
$\Uncoll_r(\calL_k,1)$ is good for all large $k$,
and the limit ray $R = \lim_{k \to \infty} \calL_k$ is irrational.
Then $W=\Uncoll_r(R,1)$ is a wonderful ray.
\end{lemma}

\begin{proof}
Using \eqref{uncoll-limit}, we see that the ray $W$ is the limit of the sequence of eventually good rays
$\Uncoll_r(\calL_k,1)$; hence $W$ is nef.  It is also irrational, since $R$ is.
\end{proof}

Our first application of this is to consider the matrix 
\[
A_n=\begin{pmatrix}
5 n + 5  &-n&  -2 n & -12 n - 12\\
5 n & 1-n& -2n& -12 n  \\
5 & -1 &-1 &-12 \\
2 & 0 & 0 & -5  
\end{pmatrix}
\]
associated to the Cremona map $\phi$ as above, 
and define integer numbers $d_{n,k}, a_{n,k}, b_{n,k}, c_{n,k}$ by
	\[
\begin{pmatrix}
d_{n,k}\\a_{n,k}\\b_{n,k}\\c_{n,k}
\end{pmatrix}=A_n^k
\begin{pmatrix}
1\\1\\0\\0
\end{pmatrix}.
\]

\begin{proposition}\label{pro:odd_good}
	For every $n\ge 0$ and every $k\ge 0$ the linear system
$$
\calP_{n,k}=\calL_{d_{n,k}}(a_{n,k},b_{n,k}^{2n},c_{n,k}^6)
$$
is a pencil of rational curves having self--intersection zero.
	For every $n\ge 2$, and every $k\ge 1$ 
the linear system
$$
\calG_{n,k}=\calL_{2d_{n,k}}(a_{n,k}^4,2b_{n,k}^{2n},2c_{n,k}^6)
$$
whose self--intersection is again zero,
is empty and all its multiples $m\calG_{n,k}$ are empty for $m\ge 1$.	
\end{proposition}

\begin{proof}
	The first claim is obvious, because $\calP_{n,k}=\phi^k(\calP_{n,0})$ 
is a Cremona transform of the pencil $\calP_{n,0}=\calL_1(1,0^{2n+6})$ 
of lines through the first point.
	
To prove that $m\calG_{n,k}$ is empty,
we collide its four points of multiplicity $ma_{n,k}$
and apply Lemma \ref{lem:uncollision}(a).
It suffices therefore to show that
$\dim \calL_{2md_n,k}(2ma_{n,k},2mb_{n,k}^{2m},2mc_{n,k}^6) < ma_{n,k}$.
This system is exactly $2m\calP_{n,k}$,
which is composed with the pencil $\calP_{n,k}$, and has dimension $2m$.
Hence it suffices to prove that $a_{n,k}>2$
(for $n \ge 2$ and $k \ge 1$).

	
	The vectors $v=(1,-1,0,-2)$ and $w=(0,-1,n,0)$ 
satisfy $vA_n=v$ and $wA_n=w$, 
so the quantities $d_{n,k}-a_{n,k}-2c_{n,k}$ and $nb_{n,k}-a_{n,k}$ 
are independent of $k$, 
and looking at the case $k=0$ we see that they equal $0$, and $-1$, respectively. 
	So to see that $a_{n,k}=1+nb_{n,k}>2$ 
it will be enough to see that $b_{n,k}>0$ for $k>1$ (because $n\ge 2$).
	On the other hand $\calP_{n,k}$ is a pencil, 
so in particular it is nef,
and therefore $d_{n,k},a_{n,k},b_{n,k},c_{n,k}\ge 0$; 
thus it only remains to prove that $b_{n,k}\ne 0$ for $k\geq1$.
	
	Again using that $a_{n,k-1}=1+nb_{n,k-1}$
and 
$d_{n,k-1}=a_{n,k-1}+2c_{n,k-1} = 1+nb_{n,k-1} + 2c_{n,k-1} $, we have

\[
\begin{pmatrix}
b_{n,k}\\c_{n,k}
\end{pmatrix}=
\begin{pmatrix}
5 & -1 &-1 &-12 \\
2 & 0 & 0 & -5  
\end{pmatrix}
\begin{pmatrix}
d_{n,k-1}\\a_{n,k-1}\\b_{n,k-1}\\c_{n,k-1}
\end{pmatrix}=
\begin{pmatrix}
4n-1 & -2 \\
2n & -1\\
\end{pmatrix}
\begin{pmatrix}
b_{n,k-1}\\c_{n,k-1}
\end{pmatrix}+
\begin{pmatrix}
4\\
2\\
\end{pmatrix}.
\]
This immediately shows that $b_{n,k}-c_{n,k}=(2n-1)b_{n,k-1}-c_{n,k-1}+2$, which allows us to prove by induction that $b_{n,k}>c_{n,k}\ge 0$. 
For $k=1$, a direct computation gives $b_{n,1}=4 > c_{n,1}=2$. 
For $k>1$, we then have $$b_{n,k}-c_{n,k}=(2n-1)b_{n,k-1}-c_{n,k-1}+2 > (2n-2)b_{n,k-1}+2>0,$$
as needed.
\end{proof}

We now apply Lemma \ref{lem:uncollwonderful} with $r=2$ and $\calL_k = \calP_{n,k}$;
we have $\calG_{n,k} = \Uncoll_2(\calP_{n,k},1)$ as rays,
and the above shows that $\calG_{n,k}$ is good for all $n \geq 2$, $k \geq 1$.
Since $\lim_{k\to\infty} \calP_{n,k}$
corresponds to the eigenvector for the dominant eigenvalue for the matrix $A_n$,
we have that this limit is $W_\mathit{odd}$.
We conclude that the limit of the $\calG_{n,k}$ is wonderful,
and this limit is equal to $\Uncoll_2(W_\mathit{odd},1)$,
which we denote by $W^+_\mathit{even}$:

\begin{corollary}\label{cor:even_wonderful}
The ray spanned by the class $\calG_{n,k}$ 
is good for every $n\ge2$, $k\ge 1$. 
Therefore its limit is wonderful; it is the ray spanned by
\begin{gather*}
W^+_\mathit{even} = 
(10 n^{2} + 8 n) H 
- n\left(3 n + 2  \alpha_n \right)\sum_{i=1}^{4}E_i \\ 
- \left( 6 n + 4 \alpha_n\right) \sum_{i=5}^{4+2n}E_i
- 2n\left(2+n-\alpha_n\right)\sum_{5+2n}^{10+2n}E_i
\end{gather*}where $\alpha_n=\sqrt{n(n-1)}$.  	
\end{corollary}

Computing the intersection with the De Fernex ray,
we see that it is negative for $n=2$; this gives us the following.

\begin{corollary}\label{sgeq14wonderfulK+}
For every even $s\ge 14$ there exist wonderful rays $R$ on $X_s$ 
with $R\cdot K_s>0$. 
For $s=14$, 
there is such a wonderful ray that is De Fernex negative.
\end{corollary}

Now to take care of the odd cases, consider the matrix 
\[
B_n=\begin{pmatrix}
8 n^{2} + 27 n + 17 & -n^{2} - 5 n - 6 & -21 n^{2} - 70 n - 42 & -2 n^{2} - 6 n \\
8 n^{2} + 19 n + 6 & -n^{2} - 4 n - 3 & -21 n^{2} - 49 n - 14 & -2 n^{2} - 4 n \\
8 n + 6 & -n - 2 & -21 n - 15 & -2 n \\
8 n + 3 & -n - 2 & -21 n - 7 & -2 n + 1	\end{pmatrix}
\]
associated to the Cremona map $\psi$ as above, 
and define integer numbers $d'_{n,k}, a'_{n,k}, b'_{n,k}, c'_{n,k}$ by
\[
\begin{pmatrix}
d'_{n,k}\\a'_{n,k}\\b'_{n,k}\\c'_{n,k}
\end{pmatrix}=B_n^k
\begin{pmatrix}
1\\1\\0\\0
\end{pmatrix}.
\]

\begin{proposition}\label{pro:even_good}
For every $n\ge 1$ and every $k\ge 0$ 
the linear system
$$
\calP'_{n,k}=\calL_{d'_{n,k}}(a'_{n,k},(b'_{n,k})^{7},(c'_{n,k})^{2n})
$$
is a pencil of rational curves, of self--intersection zero.
For every $n\ge 1$, and every $k\ge 1$
the linear system 
$$
\calG'_{n,k}=\calL_{2d'_{n,k}}((a'_{n,k})^4,(2b'_{n,k})^{7},(2c'_{n,k})^{2n})
$$
has zero self--intersection,
is empty and all its multiples $m\calG'_{n,k}$ are empty for $m\ge 1$.
\end{proposition}

\begin{proof}
The proof follows along the same lines as in the previous proposition.
The assertions about $\calP'_{n,k}$ are obvious.

To prove that $m\calG'_{n,k}$ is empty,
we collide its four points of multiplicity $ma'_{n,k}$, 
and using Lemma \ref{lem:uncollision}(a) again
we need to prove that $a'_{n,k}>2$,
whenever $n\ge 1$ and $k\ge 1$.
	
The vectors $v=(3,0,-7,-3n-2)$ and $w=(0,1,0,-n-2)$
satisfy $vB_n=v$ and $wB_n=w$,
so the quantities $3d'_{n,k}-7b'_{n,k}-(3n+2)c'_{n,k}$ and $a'_{n,k}-(n+2)c'_{n,k}$ 
are independent of $k$, 
and looking at the case $k=0$ we see that they equal $3$, and $1$, respectively. 
So to see that $a'_{n,k}=1+(n+2)c'_{n,k}>2$
it will be enough to see that $c'_{n,k}>0$ for $k>1$.
On the other hand $\calP'_{n,k}$  is nef,
and therefore $d'_{n,k},a'_{n,k},b'_{n,k},c'_{n,k}\ge 0$; 
thus it only remains to prove that $c'_{n,k}\ne 0$ for $k\geq1$.
	
Again using that $d'_{n,k-1}=(7/3)b'_{n,k-1}+(3n+2)c'_{n,k-1}+3$ 
and $a'_{n,k}=1+(n+2)c'_{n,k}$, we have
	\[
	\begin{pmatrix}
	b'_{n,k}\\c'_{n,k}
	\end{pmatrix}=
	\begin{pmatrix}
	-\frac{7}{3}n-1 & 7n^2+\frac{16}{3}n \\
	-\frac{7}{3}n & 7n^2+\frac{7}{3}n-1\\
	\end{pmatrix}
	\begin{pmatrix}
	b'_{n,k-1}\\c'_{n,k-1}
	\end{pmatrix}+
	\begin{pmatrix}
	7n+4\\
	7n+1\\
	\end{pmatrix}.
	\]
We will use this expression to show that $3c'_{n,k}>b'_{n,k}\ge 0$, 
which will finish the proof.
For $k=1$ we have $b'_{n,1}=7n+4, c'_{n,1}=7n+1$ 
satisfying the inequality (because $n\ge1$).
For $k>1$ we can argue by induction on $k$.
Indeed, it follows from the latter matrix equality that
\begin{gather*}
	3c'_{n,k}-b'_{n,k} > \left(14 n^{2} + \frac{5}{3} n - 3\right)c'_{n,k-1}
- \left(\frac{14}{3} n - 1\right)b'_{n,k-1}>\\ 
\left(\frac{14}{3} n - 1\right)(3c'_{n,k-1}-b'_{n,k-1})>0,	
	\end{gather*}and we are done.
\end{proof}

Again Lemma \ref{lem:uncollwonderful} applies,
with $r=2$ and $\calL_k = \calP'_{n,k}$;
we also have $\calG'_{n,k} = \Uncoll_2(\calP'_{n,k},1)$ as rays.
The above shows that $\calG'_{n,k}$ is good for all $n \geq 1$, $k \geq 1$,
so that since $\lim_{k\to\infty} \calP'_{n,k} = W_{\textit{even}}$,
we have that the limit of the $\calG'_{n,k}$ is wonderful
and this limit is equal to $\Uncoll_2(W_\textit{odd},1)$,
which we denote by $W^+_\mathit{odd}$:

\begin{corollary}\label{cor:odd_wonderful}
	The ray spanned by the class $\calG'_{n,k}$ is good for every $n\ge1$, $k\ge 1$. Therefore its limit is wonderful; it is the ray spanned by
\begin{gather*}
	W^+_{\textit{odd}} = 28n (8  n^{2} + 27  n + 16) H -
	7n (n + 2)(9 n + \beta_n + 6)\sum_{i=1}^{4}E_i\\
	- 2n (21 n^{2} - 3 n\beta_n + 126 n - 2 \beta_n + 84) \sum_{i=5}^{11}E_i-
	14n(9 n + \beta_n + 6)\sum_{12}^{11+2n}E_i
\end{gather*}	where $\beta_n=\sqrt{49n^2-28}$.  	
\end{corollary}

A computation shows that it is De Fernex negative for $n=1$.
Hence we have:

\begin{corollary}\label{sgeq13wonderfulK+}
	For every odd $s\ge13$ there exist wonderful rays $R$ on $X_s$ 
with $R\cdot K_s>0$. 
For $s=13$, there is  such a wonderful ray that is De Fernex negative.
\end{corollary}

\section{De Fernex negative wonderful rays}
\label{sec:deFernex_negative}

In the previous section we used uncollision of a point to four points on a wonderful ray,
to obtain a wonderful ray on $K_s^+$; 
in the initial cases of each sequence 
($n=2$ in Corollary \ref{cor:even_wonderful} 
or $n=1$ in Corollary \ref{cor:odd_wonderful}) 
the resulting ray is De Fernex negative,
but for larger $n$ the multiplicities obtained become too inhomogeneous 
and the rays become De Fernex positive. 
This can be remedied by using an uncollision to a variable number of points, 
to obtain infinite sequences of De Fernex negative wonderful rays
at the price of covering only some special values of $s$.

\begin{proposition}
	For every $n\ge 1$, and every $k\ge 1$  the linear systems
\begin{align*}
	\calG''_{n,k}&=
\calL_{(n+1)d'_{n,k}}((a'_{n,k})^{(n+1)^2},((n+1)b'_{n,k})^{7},((n+1)c'_{n,k})^{2n}) \\
	\calG'''_{n,k}&=
\calL_{(n+2)d'_{n,k}}((a'_{n,k})^{(n+2)^2},((n+2)b'_{n,k})^{7},((n+2)c'_{n,k})^{2n})
\end{align*}
have zero self--intersection, are empty and all of their multiples are empty.	
\end{proposition}
\begin{proof}
The claim for $\calG''_{1,k}=\calG'_{1,k}$ has already been proved.

For $m\calG'''_{1,k}$ and $m\calG''_{2,k}$,
we collide the $9$ points of multiplicity $ma'_{n,k}$
to a point of multiplicity $3ma'_{n,k}$ with $\alpha = 3ma'_{n,k}$ matching conditions,
as explained in Section \ref{sec:collision};
the resulting system is $3m\calP'_{n,k}$
with $\alpha$ additional conditions.
Since we already proved that $a'_{n,k}>2$, we get
$\alpha=3ma'_{n,k}>6m>\dim(3m\calP'_{n,k})$,
we can apply Lemma \ref{lem:uncollision}(b),
and conclude that the system is empty.

For all other cases,
write $r=n+1$ (in the case of $\calG''$) 
or $r=n+2$ (in the case of $\calG'''$),
and collide the $r^2$ points of multiplicity $ma'_{n,k}$
to a point of multiplicity $t>rma'_{n,k}$, because $r>3$.
The resulting linear system is the subsystem of $rm\calP'_{n,k}$
formed by the curves with a point of multiplicity $t$ at the first point.
But curves in the rational pencil $\calP'_{n,k}$
have multiplicity exactly $a'_{n,k}$ at the first point,
so curves in $rm\calP'_{n,k}$,
being sums of $rm$ curves in $\calP'_{n,k}$,
have multiplicity at most $rma_{n,k}<t$;
we conclude that collided system is empty,
and hence by Lemma \ref{lem:uncollision}(c), we have the result.
\end{proof}

Using the notation we've introduced,
we see that $\calG''_{n,k} = \Uncoll_{n+1}(\calP'_{n,k},1)$
and $\calG'''_{n,k} = \Uncoll_{n+2}(\calP'_{n,k},1)$ as rays.
The above shows that these systems are good.
The limit of the $\calP'_{n,k}$ systems is the irrational ray $W_\textit{even}$.
Hence using Lemma \ref{lem:uncollwonderful},
we conclude that the two formal uncollisions of the limit of the $W_\textit{even}$ rays
are wonderful.  We denote these two uncollisions by $W^+_{sq4}$ and $W^+_{sq2}$,
and so we have the following.

\begin{corollary}\label{cor:sparse_wonderful}
	The rays spanned by the classes $\calG''_{n,k}$ and $\calG'''_{n,k}$
are good for every $n\ge1$, $k\ge 1$.
Therefore their limits for $k\to \infty$ are wonderful; these are the rays spanned by
\begin{gather*}
	W^+_{sq4} = 14n (8  n^{2} + 27  n + 16) L -
	7n \frac{n+2}{n+1}(9 n + \beta_n + 6)\sum_{i=1}^{(n+1)^2}E_i\\
	-n (21 n^{2} - 3 n\beta_n + 126 n - 2 \beta_n + 84) \sum_{i=(n+1)^2+1}^{(n+1)^2+7}E_i-
	7n(9 n + \beta_n + 6)\sum_{(n+1)^2+8}^{(n+2)^2+4}E_i
\end{gather*}		
and
\begin{gather*}
W^+_{sq2} = 14n (8  n^{2} + 27  n + 16) L -
7n (9 n + \beta_n + 6)\sum_{i=1}^{(n+2)^2}E_i\\
-n (21 n^{2} - 3 n\beta_n + 126 n - 2 \beta_n + 84) \sum_{i=(n+2)^2+1}^{(n+2)^2+7}E_i-
7n(9 n + \beta_n + 6)\sum_{(n+2)^2+8}^{(n+3)^2+2}E_i
\end{gather*}
respectively, where $\beta_n=\sqrt{49n^2-28}$.  	
\end{corollary}

Note that the first $(n+2)^2$ and the last $2n$ points have the same multiplicity in $W^+_{sq2}$, so this class is a permutation of
	\begin{gather*}
W^+_{sq2} = 14n (8  n^{2} + 27  n + 16) L 
-n (21 n^{2} - 3 n\beta_n + 126 n - 2 \beta_n + 84) \sum_{i=1}^{7}E_i
\\- 7n(9 n + \beta_n + 6)\sum_{8}^{(n+3)^2+2}E_i.
\end{gather*}	

We were led to these examples since in the wonderful ray $W_\mathit{even}$,
the first multiplicity is exactly equal to $n+2$ times the multiplicity of the last $2n$ points.
This means that uncolliding that first point to a collection of $(n+2)^2$ points
will yield a wonderful ray with only two distinct multiplicities,
which is more uniform.
A closer examination reveals that these rays intersect the De Fernex ray negatively:

\begin{proposition}
	The classes $W^+_{sq4}$ and $W^+_{sq2}$ are De Fernex negative. 
\end{proposition}

\begin{proof}
	Since the classes $W^+_{sq4}$ and $W^+_{sq2}$ 
and the De Fernex rays are given explicitly in terms of $n$,
this is essentially a calculus exercise.
We indicate explicitly how to obtain the inequality in the case of $W^+_{sq2}$,
leaving the other class to the interested reader.
	
	The computation of the intersection product as a function of $n$
is straightforward and gives
	\begin{gather*}
	W^+_{sq2}\cdot F_{(n+3)^2+2} = 
	-63  n^{4} - 567  n^{3} - 1386  n^{2}- 756  n\\ + 14  {\left(8  n^{3} + 27  n^{2} + 16  n\right)} \sqrt{n^{2} + 6  n + 10} - 7  {\left(n^{2} + 3  n + 2\right)}\sqrt{49  n^{4} - 28  n^{2}}  .
	\end{gather*}
	We check by hand that this is negative for $n=1,2,3,4$. On the other hand it is clear that $\sqrt{n^{2} + 6  n + 10}<n+3+1/2n$ and $\sqrt{49  n^{4} - 28  n^{2}}>7n^2-3$  for all $n\ge 1$. Therefore
	\begin{gather*}
	W^+_{sq2}\cdot F_{(n+3)^2+2} <
	-63  n^{4} - 567  n^{3} - 1386  n^{2}- 756  n\\ + 14  {\left(8  n^{3} + 27  n^{2} + 16  n\right)} (n+3+1/2n) - 7  {\left(n^{2} + 3  n + 2\right)}(7n^2-3) \\ =
	7(	-7 \, n^{2} + 24 \, n + 22),
	\end{gather*}
	which is negative for $n\ge 5$.
\end{proof}

These two sets of examples give the following.

\begin{corollary}\label{dFnegativeseries}
	For every $s\ge13$ such that $s-4$ is a square and every $s\ge 18$ such that $s-2$ is a square there exist De Fernex negative wonderful rays. 
\end{corollary}

This is the final ingredient in the proof of Theorem \ref{thm:main},
which follows from Corollary \ref{sgeq10wonderful},
Corollary \ref{sgeq14wonderfulK+},
Corollary \ref{sgeq13wonderfulK+},
and Corollary \ref{dFnegativeseries}.

\end{document}